\documentclass[a4paper,10pt]{article}
\usepackage{latexsym}
\usepackage{graphicx}  
\usepackage{amsmath}   
\usepackage{array}
\usepackage{epic}
\usepackage{eepicemu}

\newtheorem{definition}{Definition}
\newtheorem{example}{Example}
\newtheorem{theorem}{Theorem}
\newtheorem{property}{Property}
\newenvironment{proof}{Proof}

\begin{document}

\begin{center}
\vspace*{13.4mm}

{\Large\bf Interval Neutrosophic Sets}\\
\vspace{9mm}
{\bf Haibin Wang\footnote{Contact Author}, Praveen Madiraju, Yanqing Zhang and Rajshekhar Sunderraman \\}
\emph{
Department of Computer Science \\
Georgia State University \\
Atlanta, Georgia 30302, USA \\
email: {\tt hwang17@student.gsu.edu}, {\tt \{cscpnmx,yzhang,raj\}@cs.gsu.edu}
}
\vspace{13.5mm}
\end{center}

\begin{abstract}
Neutrosophic set is a part of neutrosophy which studies the origin, nature,
and scope of neutralities, as well as their interactions with different
ideational spectra. Neutrosophic set is a powerful general formal
framework that has been recently proposed. However, neutrosophic set needs to
be specified from a technical point of view. To this effect,  
we define the set-theoretic operators on an instance
of neutrosophic set, we call it interval neutrosophic set (INS). We prove 
various properties of INS, which are connected to the operations
and relations over INS. Finally, we introduce and prove the convexity of 
interval
neutrosophic sets.
\end{abstract}

\vspace{9mm}

\vspace{11.7mm}

\section{Introduction}
In this section, we introduce the related works, motivation and the problems
that we are facing. 
 
\subsection{Related Works and Historical Perspective}
The concept of fuzzy sets was introduced by Zadeh in 1965~\cite{ZAD65}.
Since then fuzzy sets and fuzzy logic have been applied in many real 
applications to handle uncertainty. The traditional fuzzy set uses one real
number $\mu_A(x) \in [0, 1]$ to represent the grade of membership of fuzzy
set $A$ defined on universe $X$. Sometimes $\mu_A(x)$ itself is uncertain
and hard to be defined by a crisp value. So the concept of interval valued
fuzzy sets was proposed~\cite{TUR86} to capture the uncertainty of grade of
membership. Interval valued fuzzy set uses an interval value $[\mu_{A}^L(x),
\mu_{A}^U(x)]$ with $0 \leq \mu_{A}^L(x) \leq \mu_{A}^U(x) \leq 1$ to represent
the grade of membership of fuzzy set $A$. In some applications such as 
expert system, belief system and information fusion, we should consider not
only the truth-membership supported by the evidence but also the 
false-membership againsted by the evidence. That is beyond the scope of fuzzy
sets and interval valued fuzzy sets. In 1986, Atanassov introduced the
intuitionistic fuzzy sets~\cite{ATA86} which is a generalization of fuzzy sets
and provably equivalent to interval valued fuzzy sets. 
The intuitionistic fuzzy sets consider 
both truth-membership and false-membership. Later on, intuitionistic fuzzy
sets were extended to the interval valued intuitionistic fuzzy sets~\cite{ATA89}. The interval valued intuitionistic fuzzy set uses a pair of intervals
$[t^-, t^+], \mbox{  } 0 \leq t^-  \leq t^+ \leq 1$ and 
$[f^-, f^+], \mbox{  } 0 \leq f^- \leq f^+ \leq 1$ with $t^+ + f^+ \leq 1$ to 
describe the degree of true belief
and false belief. Because of the restriction that $t^+ + f^+ \leq 1$, 
intuitionistic fuzzy sets and interval valued intuitionistic fuzzy sets can 
only handle incomplete information not the indeterminate information and 
inconsistent information which exists commonly in belief systems. For example,
when we ask the opinion of an expert about certain statement, he or she may
say that the possibility that the statement is true is between $0.5$ and $0.7$
and the statement is false is between $0.2$ and $0.4$ and the degree that
he or she is not sure is between $0.1$ and $0.3$. Here is another example, 
suppose there are 10 votees during a voting process. In time $t_1$, three vote
``yes", two vote ``no" and five are undecided, using neutrosophic notation, it
can be expressed as $x(0.3,0.5,0.2)$. In time $t_2$, three vote ``yes", two 
vote ``no", two give up and three are undecided, it then can be expressed as
$x(0.3,0.3,02)$. That is beyond the scope of the intuitionistic fuzzy set. 
So, the notion of neutrosophic set is more general and overcomes the 
aforementioned issues. 

\subsection{Motivation}
 
In neutrosophic set, indeterminacy is quantified explicitly and \\ 
truth-membership, 
indeterminacy-membership and false-membership are independent. This assumption 
is very important in information fusion when we try to
combine the data from different sensors. 
Neutrosophy was introduced by Florentin Smarandache in 1980. ``It is a branch
of philosophy which studies the origin, nature and scope of neutralities, as
well as their interactions with different ideational spectra"~\cite{SMA99}.  
Neutrosophic set is a powerful general formal framework which generalizes the
concept of the classic set, fuzzy set~\cite{ZAD65}, 
interval valued fuzzy set~\cite{TUR86}, 
intuitionistic fuzzy set~\cite{ATA86}, 
interval valued intuitionistic fuzzy set~\cite{ATA89}, 
paraconsistent set~\cite{SMA99}, dialetheist set~\cite{SMA99}, 
paradoxist set~\cite{SMA99}, tautological set~\cite{SMA99}.
A neutrosophic set $A$ defined on universe $U$. $x = x(T, I, F) \in A$ with
$T, I$ and $F$ being the real standard or non-standard subsets of $]0^-, 1^+[$.
$T$ is the degree of true membership function in the set $A$, $I$ is the 
degree of indeterminate membership function in the set $A$ and $F$ is the
degree of false membership function in the set $A$. 

The neutrosophic set 
generalizes the above mentioned sets from philosophical point of view. 
From scientific or
engineering point of view, the neutrosophic set and set-theoretic operators
need to  be specified. Otherwise, it will be difficult to apply in the real 
applications. In this paper, we define the set-theoretic operators on an
instance of neutrosophic set called interval neutrosophic set (INS). We
call it as ``interval" because it is subclass of neutrosophic set, that is
we only consider the subunitary interval of $[0, 1]$.

\subsection{Problem Statement}

An interval neutrosophic set $A$ defined on universe $X$, $x = x(T, I, F) 
\in A$
with $T$, $I$ and $F$ being the subinterval of $[0, 1]$.
Interval neutrosophic set can represent uncertainty, imprecise, incomplete 
and inconsistent information which exist in real world. The interval 
neutrosophic set generalizes the following sets:
\begin{enumerate}
\item the \emph{classical set}, $I = \emptyset$, $\inf T = \sup T = 0$ or $1$, $\inf F = \sup F
= 0$ or $1$ and $\sup T + \sup F = 1$. 
\item the \emph{fuzzy set}, $I = \emptyset$, $\inf T = \sup T \in [0, 1]$, $\inf F = \sup F \in
[0, 1]$ and $\sup T + \sup F = 1$.
\item the \emph{interval valued fuzzy set}, $I = \emptyset$, $\inf T, \sup T, \inf F, \sup F \in
[0, 1]$, $\sup T + \inf F = 1$ and $\inf T + \sup F = 1$.
\item the \emph{intuitionistic fuzzy set}, $I = \emptyset$, $\inf T = \sup T \in [0, 1]$,
$\inf F = \sup F \in [0, 1]$ and $\sup T + \sup F \leq 1$.
\item the \emph{interval valued intuitionistic fuzzy set}, $I = \emptyset$, $\inf T, \sup T,
\inf F, \sup F \in [0, 1]$ and $\sup T + \sup F \leq 1$.
\item the \emph{paraconsistent set}, $I = \emptyset$, $\inf T = \sup T \in [0, 1]$,
$\inf F = \sup F \in [0, 1]$ and $\sup T + \sup F > 1$.
\item the \emph{interval valued paraconsistent set}, $I = \emptyset$, $\inf T, \sup T, \inf F,
\sup F \in [0, 1]$ and $\inf T + \inf F > 1$.
\end{enumerate}
  
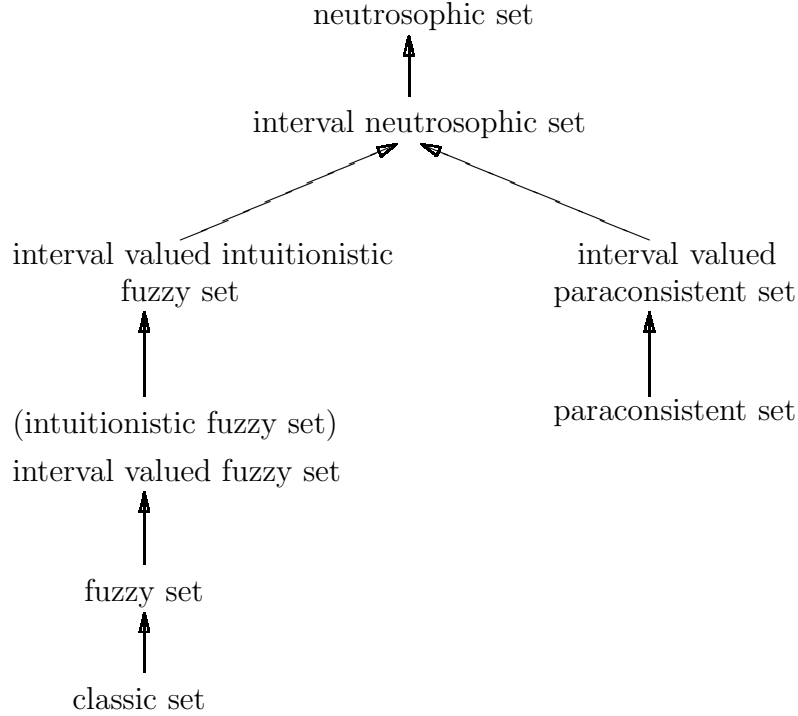
\begin{figure}[hbtp]
\centering
\setlength{\unitlength}{0.00083333in}
\begingroup\makeatletter\ifx\SetFigFont\undefined%
\gdef\SetFigFont#1#2#3#4#5{%
  \reset@font\fontsize{#1}{#2pt}%
  \fontfamily{#3}\fontseries{#4}\fontshape{#5}%
  \selectfont}%
\fi\endgroup%
{\renewcommand{\dashlinestretch}{30}
\begin{picture}(4821,4434)(0,-10)
\path(2475,3825)(2475,4200)
\blacken\path(2505.000,4080.000)(2475.000,4200.000)(2445.000,4080.000)(2505.000,4080.000)
\path(1050,2925)(2400,3525)
\blacken\path(2302.527,3448.849)(2400.000,3525.000)(2278.158,3503.678)(2302.527,3448.849)
\path(3975,2925)(2550,3525)
\blacken\path(2672.238,3506.082)(2550.000,3525.000)(2648.955,3450.784)(2672.238,3506.082)
\path(825,1950)(825,2475)
\blacken\path(855.000,2355.000)(825.000,2475.000)(795.000,2355.000)(855.000,2355.000)
\path(3975,1950)(3975,2475)
\blacken\path(4005.000,2355.000)(3975.000,2475.000)(3945.000,2355.000)(4005.000,2355.000)
\path(825,900)(825,1350)
\blacken\path(855.000,1230.000)(825.000,1350.000)(795.000,1230.000)(855.000,1230.000)
\path(825,225)(825,600)
\blacken\path(855.000,480.000)(825.000,600.000)(795.000,480.000)(855.000,480.000)
\put(1875,4275){\makebox(0,0)[lb]{\smash{{{\SetFigFont{12}{14.4}{\rmdefault}{\mddefault}{\updefault}neutrosophic set}}}}}
\put(1500,3600){\makebox(0,0)[lb]{\smash{{{\SetFigFont{12}{14.4}{\rmdefault}{\mddefault}{\updefault}interval neutrosophic set}}}}}
\put(675,2550){\makebox(0,0)[lb]{\smash{{{\SetFigFont{12}{14.4}{\rmdefault}{\mddefault}{\updefault}fuzzy set}}}}}
\put(0,2775){\makebox(0,0)[lb]{\smash{{{\SetFigFont{12}{14.4}{\rmdefault}{\mddefault}{\updefault}interval valued intuitionistic}}}}}
\put(3525,2775){\makebox(0,0)[lb]{\smash{{{\SetFigFont{12}{14.4}{\rmdefault}{\mddefault}{\updefault}interval valued}}}}}
\put(3375,2550){\makebox(0,0)[lb]{\smash{{{\SetFigFont{12}{14.4}{\rmdefault}{\mddefault}{\updefault}paraconsistent set}}}}}
\put(3375,1800){\makebox(0,0)[lb]{\smash{{{\SetFigFont{12}{14.4}{\rmdefault}{\mddefault}{\updefault}paraconsistent set}}}}}
\put(0,1425){\makebox(0,0)[lb]{\smash{{{\SetFigFont{12}{14.4}{\rmdefault}{\mddefault}{\updefault}interval valued fuzzy set}}}}}
\put(450,675){\makebox(0,0)[lb]{\smash{{{\SetFigFont{12}{14.4}{\rmdefault}{\mddefault}{\updefault}fuzzy set}}}}}
\put(375,0){\makebox(0,0)[lb]{\smash{{{\SetFigFont{12}{14.4}{\rmdefault}{\mddefault}{\updefault}classic set}}}}}
\put(0,1725){\makebox(0,0)[lb]{\smash{{{\SetFigFont{12}{14.4}{\rmdefault}{\mddefault}{\updefault}(intuitionistic fuzzy set)}}}}}
\end{picture}
}
\caption{Relationship among interval neutrosophic set and other sets}
\label{fig0}
\end{figure}

The relationship among interval neutrosophic set and other sets is 
illustrated in Fig~\ref{fig0}.

Note that $\rightarrow$ in Fig.~\ref{fig0} such as $a \rightarrow b$ means
that b is a generalization of a.

We define the set-theoretic operators on interval neutrosophic set (INS).
Various properties of INS are proved, which are connected to the operations
and relations over INS.                           

The rest of paper is organized as follows. Section~\ref{section2} gives a brief 
overview of neutrosophic set. Section~\ref{section3} gives the definition of 
interval neutrosophic set and set-theoretic operations.
Section~\ref{section4} gives some
properties of set-theoretic operations. Section~\ref{section5} gives the 
definition of convexity of interval neutrosophic sets and prove some 
properties of convexity. Section~\ref{section6} concludes the
paper. To maintain a smooth flow throughout the paper, we present the proofs
to all theorems in Appendix.  

\section{Neutrosophic Set}
\label{section2}
This section gives a brief overview of concepts of neutrosophic set defined 
in~\cite{SMA99}. Here, we use different notations to express the same meaning.
Let $S_1$ and $S_2$ be two real standard or non-standard 
subsets, then $S_1 \oplus S_2 = \{x | x=s_1 + s_2, s_1 \in S_1 \mbox{ and } s_2 \in S_2\}$, $\{1^+\} \oplus S_2 = \{x | x=1^+ + s_2, s_2 \in S_2\}$. $S_1 \ominus S_2 = \{x | x = s_1 - s_2, s_1 \in S_1 \mbox{ and } s_2 \in S_2\}$, $\{1^+\} \ominus S_2 = \{x | x = 1^+ - s_2, s_2 \in S_2\}$. $S_1 \odot S_2 = \{x | x = s_1 \cdot s_2, s_1 \in S_1 \mbox{ and } s_2 \in S_2\}$. 

\begin{definition}[Neutrosophic Set]
Let $X$ be a space of points (objects), with a generic element in $X$ denoted
by $x$. \\
A neutrosophic set $A$ in $X$ is characterized by a \emph{truth-membership function}
$T_A$, a \emph{indeterminacy-membership function} $I_A$ and a 
\emph{false-membership} function $F_A$. $T_A(x), I_A(x)$ and $F_A(x)$ are
real standard or non-standard subsets of $]0^-, 1^+[$. That is 

\begin{eqnarray}
   T_A:X & \rightarrow & ]0^-, 1^+[, \\
   I_A:X & \rightarrow & ]0^-, 1^+[, \\
   F_A:X & \rightarrow & ]0^-, 1^+[.
\end{eqnarray}

There is no restriction on the sum of $T_A(x)$, $I_A(x)$ and $F_A(x)$, so 
$0^- \leq \sup T_A(x) + \sup I_A(x) + \sup F_A(x) \leq 3^+$.
\end{definition}

\begin{definition}
The \emph{complement} of a neutrosophic set $A$ is denoted by $\bar{A}$ and 
is defined by

\begin{eqnarray}
T_{\bar{A}}(x) & = & \{1^+\} \ominus T_A(x), \\
I_{\bar{A}}(x) & = & \{1^+\} \ominus I_A(x), \\
F_{\bar{A}}(x) & = & \{1^+\} \ominus F_A(x), 
\end{eqnarray}
for all $x$ in $X$.
\end{definition}

\begin{definition}[Containment]
A neutrosophic set $A$ is \emph{contained} in the other neutrosophic set $B$, 
$A \subseteq B$,
if and only if 

\begin{eqnarray}
\inf T_A(x) \leq \inf T_B(x) & , & \sup T_A(x) \leq \sup T_B(x), \\
\inf F_A(x) \geq \inf F_B(x) & , & \sup F_A(x) \geq \sup F_B(x),
\end{eqnarray}
for all $x$ in $X$. 

\end{definition}

\begin{definition}[Union]
The \emph{union} of two neutrosophic sets $A$ and $B$ 
is a neutrosophic set
$C$, written as $C = A \cup B$, whose truth-membership, 
indeterminacy-membership and false-membership
functions are related to those of $A$ and $B$ by

\begin{eqnarray}
T_C(x) & = & T_A(x) \oplus T_B(x) \ominus T_A(x) \odot T_B(x), \\
I_C(x) & = & I_A(x) \oplus I_B(x) \ominus I_A(x) \odot I_B(x), \\
F_C(x) & = & F_A(x) \oplus F_B(x) \ominus F_A(x) \odot F_B(x),
\end{eqnarray}
for all $x$ in $X$.
\end{definition}

\begin{definition}[Intersection]
The \emph{intersection} of two neutrosophic sets $A$ and $B$ 
is a neutrosophic set $C$, written as $C = A \cap B$, whose truth-membership,
indeterminacy-membership and
false-membership functions are related to those of $A$ and $B$ by

\begin{eqnarray}
T_C(x) & = & T_A(x) \odot T_B(x), \\
I_C(x) & = & I_A(x) \odot I_B(x), \\
F_C(x) & = & F_A(x) \odot F_B(x),
\end{eqnarray}
for all $x$ in $X$.

\end{definition}

\begin{definition}[Difference]
The \emph{difference} of two neutrosophic sets $A$ and $B$ is a neutrosophic
set $C$, writeen as $C = A \setminus B$, whose truth-membership,
indeterminacy-membership and false-membership functions are related to 
those of $A$ and $B$ by

\begin{eqnarray} 
T_C(x) & = & T_A(x) \ominus T_A(x) \odot T_B(x), \\
I_C(x) & = & I_A(x) \ominus I_A(x) \odot I_B(x), \\
F_C(x) & = & F_A(x) \ominus F_A(x) \odot F_B(x),
\end{eqnarray}
for all $x$ in $X$.

\end{definition}

\begin{definition}[Cartesian Product]
Let $A$ be the neutrosophic set defined on universe $E_1$ and $B$ be the 
neutrosophic set defined on universe $E_2$. If $x(T_{A}^1,I_{A}^1,F_{A}^1) \in A$ and $y(T_{A}^2,I_{A}^2,F_{A}^2) \in B$, then the \emph{cartesian product}
of two neutrosophic sets $A$ and $B$ is defined by
\begin{equation}
(x(T_{A}^1,I_{A}^1,F_{A}^1),y(T_{A}^2,I_{A}^2,F_{A}^2)) \in A \times B
\end{equation}  
\end{definition}

\section{Interval Neutrosophic Set}
\label{section3}
In this section, we present the notion of \emph{interval neutrosophic set} 
(INS).
Interval neutrosophic set (INS) is an instance of neutrosophic set which
can be used in real scientific and engineering applications. 

\begin{definition}[Interval Neutrosophic Set]
Let $X$ be a space of points (objects), with a generic element in $X$ denoted
by $x$. \\
An interval neutrosophic set (INS) $A$ in $X$ is characterized by 
truth-membership function $T_A$, indeterminacy-membership function $I_A$ and 
false-membership function $F_A$.
For each point $x$ in $X$,
$T_A(x), I_A(x), F_A(x) \subseteq [0, 1]$.
\end{definition}

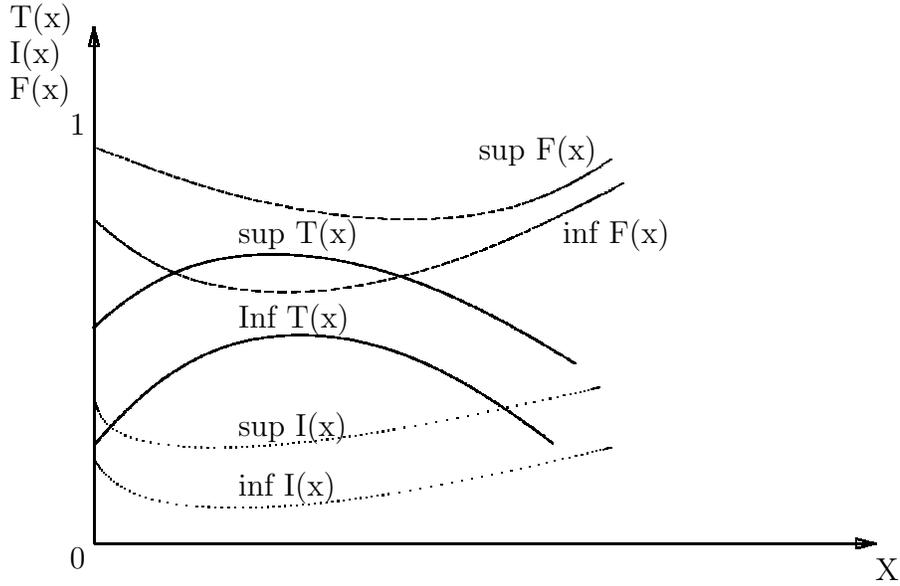
\begin{figure}[hbt]
\centering

\setlength{\unitlength}{0.00083333in}
\begingroup\makeatletter\ifx\SetFigFont\undefined%
\gdef\SetFigFont#1#2#3#4#5{%
  \reset@font\fontsize{#1}{#2pt}%
  \fontfamily{#3}\fontseries{#4}\fontshape{#5}%
  \selectfont}%
\fi\endgroup%
{\renewcommand{\dashlinestretch}{30}
\begin{picture}(5538,3609)(0,-10)
\path(525,225)(525,3450)
\blacken\path(555.000,3330.000)(525.000,3450.000)(495.000,3330.000)(555.000,3330.000)
\path(525,225)(5400,225)
\blacken\path(5280.000,195.000)(5400.000,225.000)(5280.000,255.000)(5280.000,195.000)
\path(537,852)(538,853)(539,855)
	(542,858)(547,864)(554,872)
	(564,883)(576,896)(590,912)
	(607,931)(627,952)(649,976)
	(674,1001)(701,1028)(730,1057)
	(760,1086)(793,1117)(828,1147)
	(864,1178)(901,1209)(940,1239)
	(981,1269)(1023,1297)(1068,1325)
	(1114,1352)(1162,1378)(1213,1402)
	(1267,1425)(1323,1446)(1383,1466)
	(1447,1483)(1513,1498)(1584,1510)
	(1657,1519)(1734,1525)(1812,1527)
	(1883,1525)(1955,1521)(2025,1513)
	(2094,1503)(2161,1491)(2226,1476)
	(2289,1460)(2350,1442)(2409,1423)
	(2466,1402)(2522,1380)(2576,1357)
	(2629,1333)(2680,1309)(2730,1283)
	(2780,1257)(2828,1230)(2875,1203)
	(2920,1175)(2965,1147)(3008,1120)
	(3050,1092)(3090,1066)(3129,1039)
	(3165,1014)(3199,991)(3230,968)
	(3259,948)(3285,929)(3307,912)
	(3327,898)(3343,885)(3357,875)
	(3367,867)(3375,861)(3381,857)
	(3384,854)(3386,853)(3387,852)
\path(525,1575)(526,1576)(527,1577)
	(530,1580)(535,1584)(541,1590)
	(550,1598)(561,1608)(575,1620)
	(591,1634)(610,1650)(631,1667)
	(654,1686)(680,1706)(707,1727)
	(737,1749)(768,1771)(802,1793)
	(836,1815)(873,1837)(911,1859)
	(951,1880)(993,1900)(1037,1919)
	(1083,1937)(1131,1955)(1182,1970)
	(1236,1985)(1294,1998)(1355,2009)
	(1420,2018)(1488,2025)(1561,2029)
	(1638,2031)(1718,2030)(1800,2025)
	(1872,2018)(1944,2010)(2015,1999)
	(2085,1986)(2154,1972)(2221,1956)
	(2286,1939)(2350,1921)(2411,1902)
	(2472,1882)(2530,1861)(2587,1840)
	(2643,1817)(2698,1795)(2751,1771)
	(2804,1747)(2856,1723)(2906,1699)
	(2956,1674)(3004,1649)(3052,1624)
	(3098,1600)(3142,1575)(3185,1551)
	(3227,1528)(3266,1505)(3303,1484)
	(3337,1464)(3369,1445)(3398,1428)
	(3423,1412)(3446,1399)(3466,1387)
	(3482,1377)(3495,1369)(3506,1362)
	(3514,1357)(3519,1354)(3522,1352)
	(3524,1351)(3525,1350)
\dashline{60.000}(525,2250)(526,2249)(527,2248)
	(530,2245)(535,2241)(541,2235)
	(550,2227)(561,2217)(575,2205)
	(591,2191)(610,2175)(632,2158)
	(655,2139)(681,2119)(710,2098)
	(740,2076)(772,2054)(806,2032)
	(842,2010)(880,1988)(920,1966)
	(961,1945)(1005,1925)(1051,1906)
	(1099,1888)(1151,1870)(1205,1855)
	(1263,1840)(1324,1827)(1390,1816)
	(1460,1807)(1535,1800)(1614,1796)
	(1697,1794)(1784,1795)(1875,1800)
	(1950,1806)(2026,1815)(2102,1825)
	(2176,1837)(2250,1850)(2321,1865)
	(2391,1881)(2460,1897)(2527,1915)
	(2592,1934)(2656,1953)(2718,1974)
	(2779,1994)(2839,2016)(2898,2038)
	(2956,2060)(3012,2083)(3068,2106)
	(3123,2130)(3177,2153)(3230,2177)
	(3281,2201)(3332,2224)(3381,2248)
	(3428,2271)(3473,2293)(3517,2315)
	(3558,2335)(3596,2355)(3632,2373)
	(3665,2390)(3695,2406)(3721,2420)
	(3745,2432)(3765,2443)(3781,2451)
	(3795,2459)(3806,2464)(3813,2469)
	(3819,2472)(3822,2474)(3824,2475)(3825,2475)
\dashline{60.000}(525,2700)(526,2700)(528,2699)
	(532,2697)(538,2695)(546,2692)
	(558,2687)(572,2681)(590,2674)
	(611,2666)(636,2657)(665,2646)
	(696,2634)(731,2621)(769,2607)
	(809,2592)(852,2577)(897,2560)
	(944,2544)(993,2527)(1043,2510)
	(1095,2493)(1148,2476)(1202,2459)
	(1257,2442)(1313,2426)(1370,2410)
	(1429,2394)(1488,2379)(1549,2364)
	(1611,2350)(1675,2336)(1740,2323)
	(1807,2311)(1876,2299)(1946,2289)
	(2019,2279)(2093,2270)(2169,2263)
	(2245,2257)(2323,2253)(2400,2250)
	(2489,2249)(2574,2251)(2656,2254)
	(2734,2260)(2808,2267)(2877,2276)
	(2942,2286)(3003,2297)(3060,2310)
	(3115,2323)(3166,2337)(3215,2352)
	(3261,2368)(3305,2385)(3347,2401)
	(3387,2419)(3425,2436)(3462,2454)
	(3496,2472)(3529,2489)(3560,2506)
	(3589,2523)(3616,2538)(3640,2553)
	(3662,2567)(3682,2579)(3698,2590)
	(3712,2599)(3724,2607)(3733,2613)
	(3740,2618)(3745,2621)(3748,2623)
	(3749,2624)(3750,2625)
\dottedline{45}(525,1200)(525,1199)(525,1197)
	(525,1194)(525,1190)(524,1185)
	(524,1179)(525,1171)(525,1162)
	(526,1151)(527,1140)(529,1127)
	(531,1113)(534,1099)(539,1084)
	(544,1068)(550,1052)(558,1036)
	(567,1020)(578,1004)(591,988)
	(605,972)(623,957)(642,942)
	(665,928)(690,914)(719,901)
	(752,888)(789,877)(831,866)
	(877,856)(929,847)(987,840)
	(1050,834)(1119,829)(1195,826)
	(1275,825)(1342,826)(1411,827)
	(1481,830)(1553,834)(1624,839)
	(1696,845)(1768,852)(1839,859)
	(1910,867)(1980,876)(2050,885)
	(2120,894)(2189,904)(2258,915)
	(2326,926)(2394,937)(2461,949)
	(2528,961)(2595,973)(2661,985)
	(2727,998)(2792,1010)(2857,1023)
	(2920,1036)(2983,1049)(3044,1061)
	(3104,1074)(3162,1086)(3217,1098)
	(3271,1110)(3322,1121)(3370,1131)
	(3415,1141)(3456,1150)(3494,1159)
	(3528,1167)(3558,1173)(3585,1179)
	(3607,1185)(3626,1189)(3641,1192)
	(3653,1195)(3662,1197)(3668,1198)
	(3672,1199)(3674,1200)(3675,1200)
\dottedline{45}(525,750)(526,748)(527,746)
	(529,743)(532,739)(535,733)
	(540,726)(547,717)(554,707)
	(563,696)(573,684)(585,670)
	(599,656)(614,642)(631,627)
	(649,612)(670,597)(693,582)
	(718,567)(745,553)(776,539)
	(809,526)(846,513)(887,501)
	(932,491)(982,481)(1038,472)
	(1099,464)(1166,458)(1240,453)
	(1321,450)(1408,449)(1500,450)
	(1571,452)(1644,456)(1717,460)
	(1791,466)(1865,472)(1938,479)
	(2011,487)(2082,495)(2153,504)
	(2223,514)(2292,524)(2361,534)
	(2428,545)(2495,556)(2562,568)
	(2627,580)(2692,592)(2757,604)
	(2820,617)(2883,629)(2946,642)
	(3007,655)(3068,668)(3127,681)
	(3185,694)(3242,706)(3296,718)
	(3349,730)(3398,742)(3446,753)
	(3490,763)(3531,773)(3569,781)
	(3603,789)(3633,797)(3659,803)
	(3682,808)(3701,813)(3716,817)
	(3728,820)(3737,822)(3743,823)
	(3747,824)(3749,825)(3750,825)
\put(375,75){\makebox(0,0)[lb]{\smash{{{\SetFigFont{12}{14.4}{\rmdefault}{\mddefault}{\updefault}0}}}}}
\put(5400,0){\makebox(0,0)[lb]{\smash{{{\SetFigFont{12}{14.4}{\rmdefault}{\mddefault}{\updefault}X}}}}}
\put(0,3225){\makebox(0,0)[lb]{\smash{{{\SetFigFont{12}{14.4}{\rmdefault}{\mddefault}{\updefault}I(x)}}}}}
\put(0,3450){\makebox(0,0)[lb]{\smash{{{\SetFigFont{12}{14.4}{\rmdefault}{\mddefault}{\updefault}T(x)}}}}}
\put(0,3000){\makebox(0,0)[lb]{\smash{{{\SetFigFont{12}{14.4}{\rmdefault}{\mddefault}{\updefault}F(x)}}}}}
\put(375,2775){\makebox(0,0)[lb]{\smash{{{\SetFigFont{12}{14.4}{\rmdefault}{\mddefault}{\updefault}1}}}}}
\put(1425,1575){\makebox(0,0)[lb]{\smash{{{\SetFigFont{12}{14.4}{\rmdefault}{\mddefault}{\updefault}Inf T(x)}}}}}
\put(1425,2100){\makebox(0,0)[lb]{\smash{{{\SetFigFont{12}{14.4}{\rmdefault}{\mddefault}{\updefault}sup T(x)}}}}}
\put(2925,2625){\makebox(0,0)[lb]{\smash{{{\SetFigFont{12}{14.4}{\rmdefault}{\mddefault}{\updefault}sup F(x)}}}}}
\put(3450,2100){\makebox(0,0)[lb]{\smash{{{\SetFigFont{12}{14.4}{\rmdefault}{\mddefault}{\updefault}inf F(x)}}}}}
\put(1425,900){\makebox(0,0)[lb]{\smash{{{\SetFigFont{12}{14.4}{\rmdefault}{\mddefault}{\updefault}sup I(x)}}}}}
\put(1425,525){\makebox(0,0)[lb]{\smash{{{\SetFigFont{12}{14.4}{\rmdefault}{\mddefault}{\updefault}inf I(x)}}}}}
\end{picture}
}
\caption{Illustration of interval neutrosophic set in $R^1$}
\label{fig1}
\end{figure}

An interval neutrosophic set (INS) in $R^1$ is illustrated in 
Fig.~\ref{fig1}.

When $X$ is continuous, an INS $A$ can be written as

\begin{equation}
A = \int_{X} \langle T(x), I(x), F(x) \rangle / x, \mbox{ } x \in X
\end{equation}

When $X$ is discrete, an INS $A$ can be written as

\begin{equation}
A = \sum_{i=1}^{n} \langle T(x_i), I(x_i), F(x_i) \rangle / x_i, \mbox{ } x_i \in X  
\end{equation}

Consider parameters such as capability, trusthworthiness and price of semantic
Web services. These parameters are commonly used to define quality of service
of semantic Web services.
In this section, we will use the evaluation of quality of service of semantic
Web services~\cite{WZR04} as running example to illustrate every set-theoretic
operation on interval neutrosophic set.

\begin{example}
\label{example1}
Assume that $X = [x_1, x_2, x_3]$. $x_1$ is capability, $x_2$ is trustworthiness and $x_3$ is price. The values of $x_1, x_2$ and $x_3$ are in $[0,1]$. They
are obtained from the questionaire of some domain experts, their option could
be degree of good, degree of indeterminacy and degree of poor.
$A$ is an interval  neutrosophic set 
of $X$ defined by \\
$A = \langle [0.2,0.4],[0.3,0.5],[0.3,0.5] \rangle/x_1 + \langle [0.5,0.7],[0,0.2],[0.2,0.3] \rangle/x_2 + \\ 
\langle [0.6,0.8],[0.2,0.3],[0.2,0.3] \rangle/x_3$. \\ 
$B$ is an interval neutrosophic set of $X$ defined by \\
$B = \langle [0.5,0.7],[0.1,0.3],[0.1,0.3] \rangle/x_1 + \langle [0.2,0.3],[0.2,0.4],[0.5,0.8] \rangle/x_2 + \\
\langle [0.4,0.6],[0,0.1],[0.3,0.4] \rangle/x_3$. \\
\end{example}

\begin{definition}
An interval neutrosophic set $A$ is \emph{empty} if and only if its 
$\inf T_A(x) = \sup T_A(x) = 0$, $\inf I_A(x) = \sup I_A(x) = 1$ and $\inf F_A(x) = \sup T_A(x) = 0$,
for all $x$ in $X$.
\end{definition}

We now present the set-theoretic operators on
interval neutrosophic set. 

\begin{definition}[Complement]
The \emph{complement} of an interval set $A$ is denoted by 
$\bar{A}$ and is defined by
\begin{eqnarray}
T_{\bar{A}}(x) & = & F_A(x), \\
\inf I_{\bar{A}}(x) & = & 1 - \sup I_A(x), \\
\sup I_{\bar{A}}(x) & = & 1 - \inf I_A(x), \\
F_{\bar{A}}(x) & = & T_A(x),
\end{eqnarray}
for all $x$ in $X$.
\end{definition} 

\begin{example}
\label{example2}
Let $A$ be the interval neutrosophic set defined in Example~\ref{example1}. 
Then,
${\bar{A}} = \langle [0.3,0.5],[0.5,0.7],[0.2,0.4] \rangle/x_1 + \langle [0.2,0.3],[0.8,1.0],
[0.5,0.7] \rangle/x_2 + \\
\langle [0.2,0.3],[0.7,0.8],[0.6,0.8] \rangle/x_3$. \\
 
\end{example}

\begin{definition}[Containment]
An interval neutrosophic set $A$ is \emph{contained} in the other interval 
neutrosophic set $B$, $A \subseteq B$, if and only if
\begin{eqnarray}
\inf T_A(x) \leq \inf T_B(x) & , & \sup T_A(x) \leq \sup T_B(x), \\
\inf I_A(x) \geq \inf I_B(x) & , & \sup I_A(x) \geq \sup I_B(x), \\
\inf F_A(x) \geq \inf F_B(x) & , & \sup F_A(x) \geq \sup F_B(x),
\end{eqnarray}
for all $x$ in $X$. 
\end{definition}

Note that by the definition of containment, $X$ is partial order not linear
order. For example, let $A$ and $B$ be the interval neutrosophic sets
defined in Example~\ref{example1}. Then, $A$ is not contained in $B$ and 
$B$ is not contained in $A$.
 
\begin{definition}
Two interval neutrosophic sets $A$ and $B$ are \emph{equal}, written as 
$A = B$, if and only if $A \subseteq B$ and $B \subseteq A$
\end{definition}

\begin{definition}[Union]
The \emph{union} of two interval neutrosophic sets $A$ and $B$ is an 
interval neutrosophic set $C$, written as $C = A \cup B$, whose 
truth-membership, indeterminacy-membership and false-membership functions 
are related to those
of $A$ and $B$ by
\begin{eqnarray}
\inf T_C(x) & = & \max(\inf T_{A}(x), \inf T_{B}(x)), \\
\sup T_C(x) & = & \max(\sup T_{A}(x), \sup T_{B}(x)), \\
\inf I_C(x) & = & \min(\inf I_A(x), \inf I_B(x)), \\
\sup I_C(x) & = & \min(\sup I_A(x), \sup I_B(x)), \\
\inf F_C(x) & = & \min(\inf F_A(x), \inf F_B(x)), \\
\sup F_C(x) & = & \min(\sup F_A(x), \sup F_B(x)),
\end{eqnarray}
for all $x$ in $X$. 
\end{definition}

\begin{example}
\label{example3}
Let $A$ and $B$ be the interval neutrosophic sets defined in 
Example~\ref{example1}.
Then,
$A \cup B = \langle [0.5,0.7],[0.1,0.3],[0.1,0.3] \rangle/x_1 + \\ 
\langle [0.5,0.7],
[0,0.2],
[0.2,0.3] \rangle/x_2 + 
\langle [0.6,0.8],[0,0.1],[0.2,0.3] \rangle/x_3$. \\
\end{example}

The intuition behind the union operator is that if one of elements in $A$
and $B$ is true then it is true in $A \cup B$, only both are
indeterminate and false in $A$ and $B$ then it is indeterminate and false
in $A \cup B$. The other operators should be understood similarly.

\begin{theorem}
$A \cup B$ is the smallest interval neutrosophic set containing 
both $A$ and $B$.
\end{theorem}

\begin{definition}[Intersection]
The \emph{intersection} of two interval neutrosophic sets $A$ and $B$ is 
an interval neutrosophic set $C$, written as $C = A \cap B$, 
whose truth-membership, indeterminacy-membership functions and 
false-membership functions are related to those
of $A$ and $B$ by
\begin{eqnarray}
\inf T_{C}(x) & = & \min(\inf T_{A}(x), \inf T_{B}(x)), \\
\sup T_{C}(x) & = & \min(\sup T_{A}(x), \sup T_{B}(x)), \\
\inf I_C(x)   & = & \max(\inf I_A(x), \inf I_B(x)), \\
\sup I_C(x)   & = & \max(\sup I_A(x), \sup I_B(x)), \\
\inf F_{C}(x) & = & \max(\inf F_{A}(x), \inf F_{B}(x)), \\
\sup F_{C}(x) & = & \max(\sup F_{A}(x), \sup F_{B}(x)), 
\end{eqnarray}
for all $x$ in $X$.
\end{definition}

\begin{example}
\label{example4}
Let $A$ and $B$ be the interval neutrosophic sets defined in
Example~\ref{example1}.
Then,
$A \cap B = \langle [0.2,0.4],[0.3,0.5],[0.3,0.5] \rangle/x_1 + \\ 
\langle [0.2,0.3],
[0.2,0.4],
[0.5,0.8] \rangle/x_2 + 
\langle [0.4,0.6],[0.2,0.3],[0.3,0.4] \rangle/x_3$. \\
\end{example}

\begin{theorem}
$A \cap B$ is the largest interval neutrosophic set contained in 
both $A$ and $B$.
\end{theorem}

\begin{definition}[Difference]
The \emph{difference} of two interval neutrosophic sets $A$ and $B$ 
is an interval neutrosophic 
set $C$, writeen as $C = A \setminus B$, whose truth-membership,
indeterminacy-membership and false-membership functions are related to
those of $A$ and $B$ by         
\begin{eqnarray}
\inf T_C(x) & = & \min(\inf T_A(x), \inf F_B(x)), \\
\sup T_C(x) & = & \min(\sup T_A(x), \sup F_B(x)), \\
\inf I_C(x) & = & \max(\inf I_A(x), 1 - \sup I_B(x)), \\
\sup I_C(x) & = & \max(\sup I_A(x), 1 - \inf I_B(x)), \\
\inf F_C(x) & = & \max(\inf F_A(x), \inf T_B(x)), \\
\sup F_C(x) & = & \max(\sup F_A(x), \sup T_B(x)),
\end{eqnarray}
for all $x$ in $X$.  
\end{definition}

\begin{example}
\label{example5}
Let $A$ and $B$ be the interval neutrosophic sets defined in
Example~\ref{example1}.
Then,
$A \setminus B = \langle [0.1,0.3],[0.7,0.9],[0.5,0.7] \rangle/x_1 + \\ 
\langle [0.5,0.7],
[0.6,0.8],
[0.2,0.3] \rangle/x_2 + 
\langle [0.3,0.4],[0.9,1.0],[0.4,0.6] \rangle/x_3$. \\
\end{example}

\begin{theorem}
$A \subseteq B \leftrightarrow \bar{B} \subseteq \bar{A}$
\end{theorem}

\begin{definition}[Addition]
The \emph{addition} of two interval neutrosophic sets $A$ and $B$
is an interval neutrosophic
set $C$, writeen as $C = A + B$, whose truth-membership,
indeterminacy-membership and false-membership functions are related to
those of $A$ and $B$ by
\begin{eqnarray}
\inf T_C(x) & = & \min(\inf T_A(x) + \inf T_B(x), 1), \\
\sup T_C(x) & = & \min(\sup T_A(x) + \sup T_B(x), 1), \\
\inf I_C(x) & = & \min(\inf I_A(x) + \inf I_B(x), 1), \\
\sup I_C(x) & = & \min(\sup I_A(x) + \sup I_B(x), 1), \\
\inf F_C(x) & = & \min(\inf F_A(x) + \inf F_B(x), 1), \\
\sup F_C(x) & = & \min(\sup F_A(x) + \sup F_B(x), 1),
\end{eqnarray}
for all $x$ in $X$.
\end{definition}  

\begin{example}
\label{example6}
Let $A$ and $B$ be the interval neutrosophic sets defined in
Example~\ref{example1}.
Then,
$A + B = \langle [0.7,1.0],[0.4,0.8],[0.4,0.8] \rangle/x_1 + \\
\langle [0.7,1.0],
[0.2,0.6],
[0.7,1.0] \rangle/x_2 + 
\langle [1.0,1.0],[0.2,0.4],[0.5,0.7] \rangle/x_3$. \\
\end{example}

\begin{definition}[Cartesian product]
The \emph{cartesian product} of two interval neutrosophic sets $A$ defined 
on universe $X_1$ and $B$ defined on universe $X_2$
is an interval neutrosophic
set $C$, writeen as $C = A \times B$, whose truth-membership,
indeterminacy-membership and false-membership functions are related to
those of $A$ and $B$ by
\begin{eqnarray}
\inf T_C(x,y) & = &\inf T_A(x) + \inf T_B(y) - \inf T_A(x) \cdot \inf T_B(y),\\
\sup T_C(x,y) & = &\sup T_A(x) + \sup T_B(y) - \sup T_A(x) \cdot \sup T_B(y),\\
\inf I_C(x,y) & = & \inf I_A(x) \cdot \sup I_B(y), \\
\sup I_C(x,y) & = & \sup I_A(x) \cdot \sup I_B(y), \\
\inf F_C(x,y) & = & \inf F_A(x) \cdot \inf F_B(y), \\
\sup F_C(x,y) & = & \sup F_A(x) \cdot \sup F_B(y),
\end{eqnarray}
for all $x$ in $X_1$, $y$ in $X_2$.
\end{definition}  

\begin{example}
\label{example7}
Let $A$ and $B$ be the interval neutrosophic sets defined in
Example~\ref{example1}.
Then,
$A \times B = \langle [0.6,0.82],[0.03,0.15],[0.03,0.15] \rangle/x_1 + \\
\langle [0.6,0.79],
[0,0.08],
[0.1,0.24] \rangle/x_2 + 
\langle [0.76,0.92],[0,0.03],[0.03,0.12] \rangle/x_3$. \\
\end{example}

\begin{definition}[Scalar multiplication]
The \emph{scalar multiplication} of interval neutrosophic set $A$ 
is an interval neutrosophic
set $B$, writeen as $B = a \cdot A $, whose truth-membership,
indeterminacy-membership and false-membership functions are related to
those of $A$ by
\begin{eqnarray}
\inf T_B(x) & = & \min(\inf T_A(x) \cdot a, 1), \\
\sup T_B(x) & = & \min(\sup T_A(x) \cdot a, 1), \\
\inf I_B(x) & = & \min(\inf I_A(x) \cdot a, 1), \\
\sup I_B(x) & = & \min(\sup I_A(x) \cdot a, 1), \\
\inf F_B(x) & = & \min(\inf F_A(x) \cdot a, 1), \\
\sup F_B(x) & = & \min(\sup F_A(x) \cdot a, 1),
\end{eqnarray}
for all $x$ in $X$, $a \in R^+$.
\end{definition}  

\begin{definition}[Scalar division]
The \emph{scalar division} of interval neutrosophic set $A$
is an interval neutrosophic
set $B$, writeen as $B = a \cdot A $, whose truth-membership,
indeterminacy-membership and false-membership functions are related to
those of $A$ by
\begin{eqnarray}
\inf T_B(x) & = & \min(\inf T_A(x) / a, 1), \\
\sup T_B(x) & = & \min(\sup T_A(x) / a, 1), \\
\inf I_B(x) & = & \min(\inf I_A(x) / a, 1), \\
\sup I_B(x) & = & \min(\sup I_A(x) / a, 1), \\
\inf F_B(x) & = & \min(\inf F_A(x) / a, 1), \\
\sup F_B(x) & = & \min(\sup F_A(x) / a, 1),
\end{eqnarray}
for all $x$ in $X$, $a \in R^+$.
\end{definition}             

Now we will define two operators: truth-favorite ($\triangle$) and 
false-favorite ($\nabla$) to remove the indeterminacy in the interval 
neutrosophic sets and transform it into interval valued intuitionistic
fuzzy sets or interval valued paraconsistent sets. These two operators are
unique on interval neutrosophic sets.

\begin{definition}[Truth-favorite]
The \emph{truth-favorite} of interval neutrosophic set $A$ is an interval 
neutrosophic set $B$, writeen as $B = \triangle A$, whose truth-membership
and false-membership functions are related to those of $A$ by
\begin{eqnarray}
\inf T_B(x) & = & \min(\inf T_A(x) + \inf I_A(x), 1), \\
\sup T_B(x) & = & \min(\sup T_A(x) + \sup I_A(x), 1), \\
\inf I_B(x) & = & 0, \\
\sup I_B(x) & = & 0, \\
\inf F_B(x) & = & \inf F_A(x), \\
\sup F_B(x) & = & \sup F_A(x),
\end{eqnarray}
for all $x$ in $X$.
\end{definition}

\begin{example}
\label{example8}
Let $A$ be the interval neutrosophic set defined in
Example~\ref{example1}.
Then,
$\triangle A  = \langle [0.5,0.9],[0,0],[0.3,0.5] \rangle/x_1 +
\langle [0.5,0.9],
[0,0],
[0.2,0.3] \rangle/x_2 + \\
\langle [0.8,1.0],[0,0],[0.2,0.3] \rangle/x_3$. \\
\end{example}

The purpose of truth-favorite operator is to evaluate the maximum of
degree of truth-membership of interval neutrosophic set.

\begin{definition}[False-favorite]
The \emph{false-favorite} of interval neutrosophic set $A$ is an interval 
neutrosophic set $B$, writeen as $B = \nabla A$, whose truth-membership
and false-membership functions are related to those of $A$ by
\begin{eqnarray}
\inf T_B(x) & = & \inf T_A(x), \\
\sup T_B(x) & = & \sup T_A(x), \\
\inf I_B(x) & = & 0, \\
\sup I_B(x) & = & 0, \\
\inf F_B(x) & = & \min(\inf F_A(x) + \inf I_A(x), 1), \\
\sup F_B(x) & = & \min(\sup F_A(x) + \sup I_A(x), 1),
\end{eqnarray}
for all $x$ in $X$.
\end{definition}

\begin{example}
\label{example9}
Let $A$ be the interval neutrosophic set defined in
Example~\ref{example1}.
Then,
$\nabla A  = \langle [0.2,0.4],[0,0],[0.6,1.0] \rangle/x_1 +
\langle [0.5,0.7],
[0,0],
[0.2,0.5] \rangle/x_2 + \\
\langle [0.6,0.8],[0,0],[0.4,0.6] \rangle/x_3$. \\
\end{example}

The purpose of false-favorite operator is to evaluate the maximum of
degree of false-membership of interval neutrosophic set.

\begin{theorem} 
For every two interval neutrosophic sets $A$ and $B$: \\  
\begin{enumerate}
\item $\triangle (A \cup B) \subseteq \triangle A \cup \triangle B$
\item $\triangle A \cap \triangle B \subseteq \triangle (A \cap B)$
\item $\nabla A \cup \nabla B \subseteq \nabla (A \cup B)$
\item $\nabla (A \cap B) \subseteq \nabla A \cap \nabla B$ 
\end{enumerate}
\end{theorem}

\section{Properties of Set-theoretic Operators}
\label{section4}
In this section, we will give some properties of set-theoretic operators 
defined on interval neutrosophic sets as in Section~\ref{section3}. The
proof of these properties is left for the readers.

\begin{property}[Commutativity]
$A \cup B = B \cup A$, $A \cap B = B \cap A$, $A + B = B + A$, $A \times B = B \times A$
\end{property}

\begin{property}[Associativity]
$A \cup (B \cup C) = (A \cup B) \cup C$, \\
$A \cap (B \cap C) = (A \cap B) \cap C$, \\
$A + (B + C) = (A + B) + C$, \\
$A \times (B \times C) = (A \times B) \times C$.
\end{property}

\begin{property}[Distributivity]
$A \cup (B \cap C) = (A \cup B) \cap (A \cup C)$,
$A \cap (B \cup C) = (A \cap B) \cup (A \cap C)$.
\end{property}

\begin{property}[Idempotency]
$A \cup A = A$, $A \cap A = A$, $\triangle \triangle A = \triangle A$, $\nabla \nabla A = \nabla A$.
\end{property}

\begin{property}
$A \cap \Phi = \Phi$, $A \cup X = X$, where $\inf T_{\Phi} = \sup T_{\Phi} = 0$,
$\inf I_{\Phi} = \sup I_{\Phi} = \inf F_{\Phi} = \sup F_{\Phi} = 1$ and 
$\inf T_{X} = \sup T_{X} = 1 $,
$\inf I_{X} = \sup I_{X} = \inf F_{X} = \sup F_{X} = 0$. 
\end{property}

\begin{property}
$\triangle (A + B) = \triangle A + \triangle B$, $\nabla (A + B) = \nabla A + \nabla B$.
\end{property}

\begin{property}
$A \cup \Psi = A$, $A \cap X = A$, where $\inf T_{\Phi} = \sup T_{\Phi} = 0 $,
$\inf I_{\Phi} = \sup I_{\Phi} = \inf F_{\Phi} = \sup F_{\Phi} = 1$ and 
$\inf T_{X} = \sup T_{X} = 1$,
$\inf I_{X} = \sup I_{X} = \inf F_{X} = \sup F_{X} = 0$.
\end{property}

\begin{property}[Absorption]
$A \cup (A \cap B) = A$, $A \cap (A \cup B) = A$
\end{property}

\begin{property}[DeMorgan's Laws]
$\overline{A \cup B} = \bar{A} \cap \bar{B}$, 
$\overline{A \cap B} = \bar{A} \cup \bar{B}$.
\end{property}

\begin{property}[Involution]
$\overline{{\overline{A}}} = A$
\end{property}

Here, we notice that by the definitions of complement, union and intersection
of interval neutrosophic set, interval neutrosophic set satisfies the
most properties of class set, fuzzy set and intuitionistic fuzzy set. Same as
fuzzy set and intuitionistic fuzzy set, it does not satisfy the principle of
middle exclude.
 
\section{Convexity of Interval Neutrosophic Set}
\label{section5}
We assume that $X$ is a real Euclidean space $E^n$ for correctness.

\begin{definition}[Convexity]
An interval neutrosophic set $A$ is convex if and only if 
\begin{eqnarray}
\inf T_A(\lambda x_1 + (1 - \lambda) x_2) & \geq & \min(\inf T_A(x_1), \inf T_A(x_2)), \\
\sup T_A(\lambda x_1 + (1 - \lambda) x_2) & \geq & \min(\sup T_A(x_1), \sup T_A(
x_2)), \\
\inf I_A(\lambda x_1 + (1 - \lambda) x_2) & \leq & \max(\inf I_A(x_1), \inf I_A(x_2)), \\
\sup I_A(\lambda x_1 + (1 - \lambda) x_2) & \leq & \max(\sup I_A(x_1), \sup I_A(
x_2)), \\
\inf F_A(\lambda x_1 + (1 - \lambda) x_2) & \leq & \max(\inf F_A(x_1), \inf F_A(x_2)), \\
\sup F_A(\lambda x_1 + (1 - \lambda) x_2) & \leq & \max(\sup F_A(x_1), \sup F_A(
x_2)),
\end{eqnarray}
for all $x_1$ and $x_2$ in $X$ and all $\lambda$ in $[0, 1]$. 
\end{definition}
Fig.~\ref{fig1}
is an illustration of convex interval neutrosophic set.

\begin{theorem}
If $A$ and $B$ are convex, so is their intersection. 
\end{theorem}

\begin{definition}[Strongly Convex]
An interval neutrosophic set $A$ is \\ strongly  convex if for any two distinct
points $x_1$ and $x_2$, and any $\lambda$ in the \\ open interval $(0, 1)$, 
\begin{eqnarray}
\inf T_A(\lambda x_1 + (1 - \lambda) x_2) & > & \min(\inf T_A(x_1), \inf T_A(
x_2)), \\
\sup T_A(\lambda x_1 + (1 - \lambda) x_2) & > & \min(\sup T_A(x_1), \sup T_A(
x_2)), \\
\inf I_A(\lambda x_1 + (1 - \lambda) x_2) & < & \max(\inf I_A(x_1), \inf I_A(
x_2)), \\
\sup I_A(\lambda x_1 + (1 - \lambda) x_2) & < & \max(\sup I_A(x_1), \sup I_A(
x_2)), \\
\inf F_A(\lambda x_1 + (1 - \lambda) x_2) & < & \max(\inf F_A(x_1), \inf F_A(
x_2)), \\
\sup F_A(\lambda x_1 + (1 - \lambda) x_2) & < & \max(\sup F_A(x_1), \sup F_A(
x_2)),
\end{eqnarray}
for all $x_1$ and $x_2$ in $X$ and all $\lambda$ in $[0, 1]$. 
\end{definition}

\begin{theorem}
If $A$ and $B$ are strongly convex, so is their instersection.
\end{theorem}

\section{Conclusions and Future Works}
\label{section6}
In this paper, we have presented an instance of neutrosophic set called
interval neutrosophic set (INS). The interval neutrosophic set is 
a generalization of classic set, fuzzy set, interval valued fuzzy set,
intuitionistic fuzzy sets, interval valued intuitionistic fuzzy set,
interval type-2 fuzzy 
set~\cite{LM00} and paraconsistent set.
The notions of inclusion, union, intersection, complement, relation, and 
composition have been defined on interval neutrosophic set. Various 
properties of set-theoretic operators have been proved. 
In the future, we will create the logic inference system based
on interval neutrosophic set and apply the theory to solve practical
applications in areas such as such as
expert system, data mining, question-answering system, bioinformatics and
database, etc.

\section*{Appendix}
\setcounter{theorem}{0}

\begin{theorem}
\label{theorem0}
$A \cup B$ is the smallest interval neutrosophic set containing
both $A$ and $B$.
\end{theorem}
\begin{proof}
Let $C = A \cup B$. $\inf T_C = \max(\inf T_{A}, \inf T_{B})$, $\inf T_C \geq
\inf T_A$, $\inf T_{C} \geq \inf T_B$.
$\sup T_{C} = \max(\sup T_A, \sup T_B$,
$\sup T_C \geq \sup T_A$, $\sup T_C \geq \sup T_B$.
$\inf I_C = \min(\inf I_A, \inf I_B)$, $\inf I_C \leq \inf I_A$,
$\inf I_C \leq \inf I_B$, \\
$\sup I_C = \min(\sup I_A, \sup I_B)$, $\sup I_C \leq \sup I_A$,
$\sup I_C \leq \sup I_B$,
$\inf F_C = \min(\inf F_A, \inf F_B)$, $\inf F_{C} \leq \inf F_{A}$,
$\inf F_C \leq \inf F_B$. \\
$\sup F_{C} = \min(\sup F_A, \sup F_B)$, $\sup F_C \leq \sup F_A$,
$\sup F_C \leq  \sup F_B$. 
That means $C$ contains both $A$ and $B$. \\
Furthermore, if $D$ is any extended vague set containing both $A$ and $B$, then
$\inf T_{D} \geq \inf T_{A}$, $\inf T_{D} \geq \inf T_B$, so
$\inf T_{D} \geq \max(\inf T_A, \inf T_B) = \inf T_{C}$.
$\sup T_{D} \geq \sup T_{A}$,
$\sup T_{D} \geq \sup T_{B}$,
so $\sup T_{D} \geq \max(\sup T_{A}, \sup T_{B}) = \sup T_{C}$.
$\inf I_D \leq \inf I_A$, $\inf I_D \leq \inf I_B$, so
$\inf I_{D} \leq \min(\inf I_A, \inf I_B) = \inf I_C$.
$\sup I_D \leq \sup I_A$, $\sup I_D \leq \sup I_B$, so
$\sup I_D \leq \min(\sup I_A, \sup I_B) = \sup I_C$.
$\inf F_{D} \leq \inf F_{A}$,
$\inf F_{D} \leq \inf F_B$, so $\inf F_{D} \leq \min(\inf F_{A}, \inf F_{B})
= \inf F_{C}$.
$\sup F_{D} \leq \sup F_{A}$, $\sup F_{D} \leq \sup F_{B}$, so
$\sup F_{D} \leq \min(\sup F_{A}, \sup F_{B}) = \sup F_{C}$.
That implies $C \subseteq D$.
\end{proof}

\begin{theorem}
$A \cap B$ is the largest interval neutrosophic set contained in
both $A$ and $B$.
\end{theorem}
\begin{proof}
The proof is analogous to the proof of theorem~\ref{theorem0}.
\end{proof}

\begin{theorem}
$A \subseteq B \leftrightarrow \bar{B} \subseteq \bar{A}$
\end{theorem}
\begin{proof}
$A \subseteq B \Leftrightarrow \inf T_A \leq \inf T_B, \sup T_A \leq \sup T_B,
\inf I_A \geq \inf I_B, \sup I_A \geq \sup I_B, \inf F_A \geq \inf F_B,
\sup F_A \geq \sup F_B
 \Leftrightarrow
\inf F_{B} \leq \inf F_{A}, \sup F_{B} \leq \sup F_{A}, \\ 1 - \sup I_B \geq 1 -
 \sup I_A, 1 - \inf I_B \geq 1 - \inf I_A, \inf T_{B} \geq \inf T_{A},
\sup T_{B} \geq \sup T_{A} \Leftrightarrow \bar{B} \subseteq \bar{A}$.
\end{proof}

\begin{theorem}
For every two interval neutrosophic sets $A$ and $B$: \\
\begin{enumerate}
\item $\triangle (A \cup B) \subseteq \triangle A \cup \triangle B$
\item $\triangle A \cap \triangle B \subseteq \triangle (A \cap B)$
\item $\nabla A \cup \nabla B \subseteq \nabla (A \cup B)$
\item $\nabla (A \cap B) \subseteq \nabla A \cap \nabla B$
\end{enumerate}
\end{theorem}
\begin{proof}
We now prove the first identity. Let $C = A \cup B$.\\
$\inf T_C(x) = \max(\inf T_A(x), \inf T_B(x))$, \\
$\sup T_C(x) = \max(\sup T_A(x), \sup T_B(x))$, \\
$\inf I_C(x) = \min(\inf I_A(x), \inf I_B(x))$, \\
$\sup I_C(x) = \min(\sup I_A(x), \sup I_B(x))$, \\
$\inf F_C(x) = \min(\inf F_A(x), \inf F_B(x))$, \\
$\sup F_C(x) = \min(\sup F_A(x), \sup F_B(x))$. \\
$\inf T_{\triangle C}(x) = \min(\inf T_C(x) + \inf I_C(x), 1)$, \\
$\sup T_{\triangle C}(x) = \min(\sup T_C(x) + \sup T_C(x), 1)$, \\
$\inf I_{\triangle C}(x) = \sup I_{\triangle C}(x) = 0$, \\
$\inf F_{\triangle C}(x) = \inf I_C(x)$, \\
$\sup F_{\triangle C}(x) = \sup I_C(x)$. \\
$\inf T_{\triangle A}(x) = \min(\inf T_A(x) + \inf I_A(x), 1)$, \\
$\sup T_{\triangle A}(x) = \min(\sup T_A{x} + \sup I_A(x), 1)$, \\
$\inf I_{\triangle A}(x) = \sup I_{\triangle A}(x) = 0$, \\
$\inf F_{\triangle A}(x) = \inf I_A(x)$, \\
$\sup F_{\triangle A}(x) = \sup I_A(x)$. \\
$\inf T_{\triangle B}(x) = \min(\inf T_B(x) + \inf I_B(x), 1)$, \\
$\sup T_{\triangle B}(x) = \min(\sup T_B{x} + \sup I_B(x), 1)$, \\
$\inf I_{\triangle B}(x) = \sup I_{\triangle B}(x) = 0$, \\
$\inf F_{\triangle B}(x) = \inf I_B(x)$, \\
$\sup F_{\triangle B}(x) = \sup I_B(x)$. \\
$\inf T_{\triangle A \cup \triangle B}(x) = \max(\inf T_{\triangle A}(x), \inf T_{\triangle B}(x))$, \\
$\sup T_{\triangle A \cup \triangle B}(x) = \max(\sup T_{\triangle A}(x), \sup T
_{\triangle B}(x))$, \\
$\inf I_{\triangle A \cup \triangle B}(x) = \sup I_{\triangle A \cup \triangle B}(x) = 0$, \\
$\inf F_{\triangle A \cup \triangle B}(x) = \min(\inf F_{\triangle A}(x), \inf F_{\triangle B}(x))$, \\ 
$\sup F_{\triangle A \cup \triangle B}(x) = \min(\inf F_{\triangle A}(x), \inf F
_{\triangle B}(x))$. \\
Because, \\
$\inf T_{\triangle (A \cup B)} \leq \inf T_{\triangle A \cup \triangle B}$, \\
$\sup T_{\triangle (A \cup B)} \leq \sup T_{\triangle A \cup \triangle B}$, \\
$\inf I_{\triangle (A \cup B)} =  \inf T_{\triangle A \cup \triangle B} = 0$,\\
$\sup I_{\triangle (A \cup B)} =  \sup T_{\triangle A \cup \triangle B} = 0$,\\
$\inf F_{\triangle (A \cup B)} =  \inf F_{\triangle A \cup \triangle B}$,\\
$\sup F_{\triangle (A \cup B)} =  \sup T_{\triangle A \cup \triangle B}$,\\
so, $\triangle (A \cup B) \subseteq \triangle A \cup \triangle B$.
The other identities can be proved in a similar manner.
\end{proof}

\begin{theorem}
\label{theorem4}
If $A$ and $B$ are convex, so is their intersection.
\end{theorem}
\begin{proof}
Let $C = A \cap B$, then \\
$\inf T_C(\lambda x_1 + (1 - \lambda) x_2) \geq \min(\inf T_A(\lambda x_1 + (1 -
 \lambda) x_2), \inf T_B(\lambda x_1 + (1 - \lambda) x_2))$,
$\sup T_C(\lambda x_1 + (1 - \lambda) x_2) \geq \min(\sup T_A(\lambda x_1 + (1 -
 \lambda) x_2), \sup T_B(\lambda x_1 + (1 - \lambda) x_2))$,
$\inf I_C(\lambda x_1 + (1 - \lambda) x_2) \leq \max(\inf I_A(\lambda x_1 + (1 -
 \lambda) x_2), \inf I_B(\lambda x_1 + (1 - \lambda) x_2))$,
$\sup I_C(\lambda x_1 + (1 - \lambda) x_2) \leq \max(\sup I_A(\lambda x_1 + (1 -
 \lambda) x_2), \sup I_B(\lambda x_1 + (1 - \lambda) x_2))$,
$\inf F_C(\lambda x_1 + (1 - \lambda) x_2) \leq \max(\inf F_A(\lambda x_1 + (1 -
 \lambda) x_2), \inf F_B(\lambda x_1 + (1 - \lambda) x_2))$,
$\sup F_C(\lambda x_1 + (1 - \lambda) x_2) \leq \max(\inf F_A(\lambda x_1 + (1 -
 \lambda) x_2), \inf F_B(\lambda x_1 + (1 - \lambda) x_2))$,
Since $A$ and $B$ are convex:
$\inf T_A(\lambda x_1 + (1 - \lambda) x_2) \geq \min(\inf T_A(x1), \inf T_A(x2))
$,
$\sup T_A(\lambda x_1 + (1 - \lambda) x_2) \geq \min(\sup T_A(x1), \sup T_A(x2))
$,
$\inf I_A(\lambda x_1 + (1 - \lambda) x_2) \leq \max(\inf I_A(x1), \inf I_A(x2))
$, \
$\sup I_A(\lambda x_1 + (1 - \lambda) x_2) \leq \\ \max(\sup I_A(x1), \sup I_A(x
2))
$,
$\inf F_A(\lambda x_1 + (1 - \lambda) x_2) \leq \\ \max(\inf F_A(x1), \inf F_A(x
2))
$,
$\sup F_A(\lambda x_1 + (1 - \lambda) x_2) \leq \\ \max(\sup F_A(x1), \sup F_A(x
2))
$, \\
$\inf T_B(\lambda x_1 + (1 - \lambda) x_2) \geq \min(\inf T_B(x1), \inf T_A(x2))
$,
$\sup T_B(\lambda x_1 + (1 - \lambda) x_2) \geq \min(\sup T_B(x1), \sup T_A(x2))
$,
$\inf I_B(\lambda x_1 + (1 - \lambda) x_2) \leq \\ \max(\inf I_B(x1), \inf I_A(x
2))
$,
$\sup I_B(\lambda x_1 + (1 - \lambda) x_2) \leq \\ \max(\sup I_B(x1), \sup I_A(x
2))
$, \\
$\inf F_B(\lambda x_1 + (1 - \lambda) x_2) \leq \max(\inf F_B(x1), \inf F_A(x2))
$, \\
$\sup F_B(\lambda x_1 + (1 - \lambda) x_2) \leq \max(\sup F_B(x1), \sup F_A(x2))
$, \\
Hence, \\
$\inf T_C(\lambda x_1 + (1 - \lambda) x_2) \geq \min(\min(\inf T_A(x_1), \inf T_
A(x_2)) \\
, \min(\inf T_B(x_1), \inf T_B(x_2))) = \min(\min(\inf T_A(x_1), \inf T_B(x_1)),
\\ \min(\inf T_A(x_2), \inf T_B(x_2))) =  \min(\inf T_C(x_1), \inf T_C(x_2))$,
$\sup T_C(\lambda x_1 + (1 - \lambda) x_2) \geq \min(\min(\sup T_A(x_1), \sup T_
A(x_2))
, \\ \min(\sup T_B(x_1), \sup T_B(x_2))) = \min(\min(\sup T_A(x_1), \sup T_B(x_1
)), \\
\min(\sup T_A(x_2), \sup T_B(x_2))) =  \min(\sup T_C(x_1), \sup T_C(x_2))$, \\
$\inf I_C(\lambda x_1 + (1 - \lambda) x_2) \leq \max(\max(\inf I_A(x_1), \\ \inf
 I_
A(x_2))
, \max(\inf I_B(x_1), \inf I_B(x_2))) = \max(\max(\inf I_A(x_1), \\ \inf I_B(x_1
)),
\max(\inf I_A(x_2), \inf I_B(x_2))) = \max(\inf I_C(x_1), \inf I_C(x_2))$, \\
$\sup I_C(\lambda x_1 + (1 - \lambda) x_2) \leq \max(\max(\sup I_A(x_1), \sup I_
A(x_2)) \\
, \max(\sup I_B(x_1), \sup I_B(x_2))) = \max(\max(\sup I_A(x_1), \sup I_B(x_1)),
\\
\max(\sup I_A(x_2), \sup I_B(x_2))) = \max(\sup I_C(x_1), \sup I_C(x_2))$, \\
$\inf F_C(\lambda x_1 + (1 - \lambda) x_2) \leq \max(\max(\inf F_A(x_1), \inf F_
A(x_2))
, \\ \max(\inf F_B(x_1), \inf F_B(x_2))) = \max(\max(\inf F_A(x_1), \inf F_B(x_1
)), \\
\max(\inf F_A(x_2), \inf F_B(x_2))) = \\ \max(\inf F_C(x_1), \inf F_C(x_2))$,
$\sup F_C(\lambda x_1 + (1 - \lambda) x_2) \\ \leq \max(\max(\sup F_A(x_1), \sup
 F_
A(x_2))
, \\ \max(\sup F_B(x_1), \sup F_B(x_2))) = \max(\max(\sup F_A(x_1), \sup F_B(x_1
)), \\
\max(\sup F_A(x_2), \sup F_B(x_2))) = \max(\sup F_C(x_1), \sup F_C(x_2))$.
\end{proof}

\begin{theorem}
If $A$ and $B$ are strongly convex, so is their instersection.
\end{theorem}
\begin{proof}
The proof is analogus to the proof of Theorem~\ref{theorem4}.
\end{proof}

\section*{Acknowledgments}
The authors would like to thank Dr. Florentin Smarandache for his valuable
suggestions.

\end{document}